\documentclass[12pt]{article}
\usepackage{amsmath,amsthm,amsfonts,a4wide}

\newtheorem{theorem}{Theorem}
\newtheorem{lemma}[theorem]{Lemma}
\newtheorem{remark}[theorem]{Remark}
\newtheorem{corollary}[theorem]{Corollary}
\let\eps\varepsilon
\newcommand{\R}{{\mathbb{R}}}
\newcommand{\N}{{\mathbb{N}}}
\newcommand{\Hh}{{\cal H}}
\newcommand{\eqn}[1]{(\ref{#1})}

\begin{document}

\title{A nonlinear fourth-order parabolic equation
and related logarithmic Sobolev inequalities\thanks{The authors
acknowledge partial
support from the Project ``Hyperbolic and Kinetic Equations'' of
the European Union, grant HPRN-CT-2002-00282, and from the
DAAD-Procope Program. The last author has been supported by the
Deutsche Forschungsgemeinschaft, grants JU359/3 (Gerhard-Hess
Award) and JU359/5 (Priority Program ``Multi-scale Problems'').}}
\author{Jean Dolbeault$^\dag$, Ivan Gentil\thanks{Ceremade (UMR CNRS 7534),
Universit\'e Paris IX-Dauphine, Place de la Lattre de Tassigny,
75775 Paris, C\'edex 16, France; e-mail:
{\tt $\{$dolbeaul,gentil$\}$@ceremade.dauphine.fr.}}, and
Ansgar J\"ungel\thanks{
Fachbereich Mathematik und Informatik, Universit{\"a}t Mainz,
Staudingerweg 9, 55099 Mainz, Germany; e-mail: {\tt
juengel@mathematik.uni-mainz.de.}}}

\date{June 16, 2004}

\maketitle

\begin{abstract}\noindent A nonlinear fourth-order parabolic equation in one space dimension with periodic boundary conditions is studied. This equation arises in the context of fluctuations of a stationary nonequilibrium interface and in the modeling of quantum semiconductor devices. The existence of global-in-time non-negative weak solutions is shown. A criterion for the uniqueness of non-negative weak solutions is given. Finally, it is proved that the solution converges exponentially fast to its mean value in the ``entropy norm'' using a new optimal logarithmic Sobolev inequality for higher derivatives.
\end{abstract}

{\small \bigskip\noindent {\bf AMS Classification.} 35K35, 35K55, 35B40.

}


\section{Introduction}

This paper is concerned with the study of some properties of weak
solutions to a nonlinear fourth-order equation with periodic
boundary conditions and related logarithmic Sobolev inequalities.
More precisely, we consider the problem
\begin{equation}
  u_t + (u(\log u)_{xx})_{xx} = 0, \quad u(\cdot,0)=u_0\ge 0 \quad\mbox{in }
    S^1, \label{equ}
\end{equation}
where $S^1$ is the one-dimensional torus parametrized by a variable
$x\in[0,L]$.

Recently equation \eqref{equ} has attracted the interest of many
mathematicians since it possesses some remarkable properties, e.g.,
the solutions are non-negative and there are several Lyapunov functionals.
For instance, a formal calculation shows that the
{\em entropy} is non-increasing:
\begin{equation}
  \frac{d}{dt}\int_{S^1}u(\log u-1)dx + \int_{S^1}u\left|(\log u)_{xx}\right|^2 dx = 0.
    \label{en1}
\end{equation}
Another example of a Lyapunov functional is $\int_{S^1}(u-\log u)dx$ which formally yields
\begin{equation}
  \frac{d}{dt}\int_{S^1}(u-\log u)dx + \int_{S^1}\left|(\log u)_{xx}\right|^2dx = 0.
    \label{en2}
\end{equation}
This last estimate is used to prove that solutions
to \eqref{equ} are non-negative. Indeed, a Poincar\'e inequality shows that $\log u$ is bounded in $H^2(S^1)$
and hence in $L^\infty(S^1)$, which implies that $u\ge 0$ in
$S^1\times(0,\infty)$. We prove this result rigorously in section
\ref{sec-ex}.  Notice that the equation is of
higher order and no maximum principle argument can be employed.
For more comments on Lyapunov functionals of \eqref{equ} we refer to \cite{BLS94,CCT04}.

Equation \eqref{equ} has been first derived in the context
of fluctuations of a stationary non-equilibrium interface \cite{DLSS91}.
It also appears as a zero-temperature zero-field approximation of
the so-called quantum drift-diffusion model for semiconductors
\cite{Anc87} which can be derived by a quantum moment method from
a Wigner-BGK equation \cite{DMR04}. The first analytical result
has been presented in \cite{BLS94}; there the existence of local-in-time
classical solutions with periodic boundary conditions has been proved.
A global-in-time existence result with homogeneous Dirichlet-Neumann
boundary conditions has been obtained in \cite{JuPi00}. However, up to now,
no global-in-time existence result is available for the problem
\eqref{equ}.

The long-time behavior of solutions has been studied in
\cite{CCT04} using periodic boundary conditions, in \cite{JuTo03}
with homogeneous Dirichlet-Neumann boundary conditions and finally,
in \cite{GJT04} employing non-homogeneous Dirichlet-Neumann
boundary conditions. In particular, it has been shown that the
solutions converge exponentially fast to their steady state. The decay rate
has been numerically computed in \cite{CJT03}. We also mention the
work \cite{JuPi01} in which a positivity-preserving numerical scheme
for the quantum drift-diffusion model has been proposed.

In the last years the question of non-negative or
positive solutions of fourth-order parabolic equations has also been
investigated in the context of lubrication-type equations, like the
thin film equation
$$
  u_t + (f(u)u_{xxx})_x = 0
$$
(see, e.g., \cite{BeFr90,Ber98}), where typically, $f(u)=u^\alpha$ for some
$\alpha>0$. This equation is of degenerate type which makes the analysis
easier than for \eqref{equ}, at least concerning the positivity property.

\medskip In this paper we show the following results. First, the existence of
global-in-time weak solutions is shown under a rather weak condition
on the initial datum $u_0$. We only assume that $u_0\ge 0$ is measurable
and such that $\int_{S^1}(u_0-\log u_0)dx<\infty$. Compared to \cite{BLS94},
we do not impose any smallness condition on $u_0$. We are able to prove that
the solution is non-negative. The main idea of the proof consists in performing an
exponential change of unknowns of the form $u=e^y$ and to solve a
semi-discrete approximate problem. An estimate similar to \eqref{en2}
and a Poincar\'e inquality provide $H^2$ bounds for $\log u=y$, which are uniform
in the approximation parameter. Performing the limit in this parameter
yields a {\em non-negative} solution to \eqref{equ}. These ideas
have been already employed in \cite{JuPi00} but here we need an additional
regularization procedure in the linearized problem in order to replace the usual Poincar\'e inequality in $H_0^1$ (see the proof of Theorem~\ref{existence} for details).

Our second result is concerned with uniqueness issues. If $u_1$ and $u_2$ are two
non-negative solutions to~\eqref{equ} satisfying some regularity
assumptions (see Theorem \ref{unique}) then $u_1=u_2$. A uniqueness
result has already been obtained in \cite{BLS94} in the class of
mild {\em positive} solutions; however, our result
allows for all {\em non-negative} solutions satisfying only a few additional
assumptions.

The third result is the exponential time decay of the solutions, i.e.,
we show that the solution constructed in Theorem \ref{existence}
converges exponentially fast to its mean value $\bar u=\int u(x,t)dx/L$:
\begin{equation}
  \int_{S^1}u(x,t)\log\left(\frac{u(x,t)}{\bar u}\right)dx \le e^{-Mt}
    \int_{S^1}u_0\log\left(\frac{u_0}{\bar u}\right)dx \quad \forall\;t>0, \label{decay}
\end{equation}
where $M=32\pi^4/L^4$. The same constant has been obtained in \cite{CCT04}
(even in the $H^1$ norm); however, our proof is based on the entropy--entropy production method and therefore much simpler. For this, we show that the entropy production term
$\int u\left|(\log u)_{xx}\right|^2dx$ in \eqref{en1} can be bounded from below
by the entropy itself yielding
$$
  \frac{d}{dt}\int_{S^1}u\log\left(\frac{u}{\bar u}\right)dx + M\int_{S^1}
    u\log\left(\frac{u}{\bar u}\right)dx \le 0.
$$
Then Gronwall's inequality gives \eqref{decay}. This argument is
formal since we only have weak solutions; we refer to Theorem
\ref{longtime} for details of the rigorous proof.

The lower bound for the entropy production is obtained through
a logarithmic Sobolev inequality in $S^1$. We show (see Theorem
\ref{logsob}) that any function $u\in H^n(S^1)$ ($n\in\N$) satisfies
\begin{equation}
  \int_{S^1} u^2\log\left(\frac{u^2}{\|u\|_{L^2(S^1)}^2}\right)dx
   \le 2\left(\frac{L}{2\pi}\right)^{2n}\int_{S^1}\left|u^{(n)}\right|^2 dx,  \label{opt}
\end{equation}
where $\|u\|_{L^2(S^1)}^2=\int u^2 dx/L$, and the constant is {\em optimal}.
As already mentioned in the case $n=2$, the proof of this result uses the
entropy--entropy production method.

\medskip The paper is organized as follows. In section \ref{sec-ex} the existence
of solutions is proved. Section \ref{sec-un} is concerned with the uniqueness
result. Then section \ref{sec-so} is devoted to the proof of the
optimal logarithmic Sobolev inequality \eqref{opt}.
Finally, in section \ref{sec-lo}, the exponential time decay
\eqref{decay} is shown.


\section{Existence of solutions}\label{sec-ex}

\begin{theorem}\label{existence}
Let $u_0:S^1\to{\mathbb R}$ be a nonnegative measurable function such that
$\int_{S^1}(u_0-\log u_0)dx<\infty$.
Then there exists a global weak solution $u$ of \eqref{equ} satisfying
\begin{eqnarray*}
  u\in L^{5/2}_{\rm loc}(0,\infty;W^{1,1}(S^1))\cap W^{1,1}_{\rm loc}(0,\infty;
    H^{-2}(S^1)), \\
  u\ge 0\quad\mbox{in }S^1\times(0,\infty), \quad
    \log u\in L^2_{\rm loc}(0,\infty;H^2(S^1)),
\end{eqnarray*}
and for all $T>0$ and all smooth test functions $\phi$,
$$
  \int_0^T\langle u_t,\phi\rangle_{H^{-2},H^2}dt + \int_0^T\int_{S^1}
    u(\log u)_{xx}\phi_{xx} dx dt = 0.
$$
The initial datum is satisfied in the sense of $H^{-2}(S^1):=
(H^2(S^1))^*$.
\end{theorem}

\proof We first transform \eqref{equ} by introducing the new variable
$u=e^y$ as in \cite{JuPi00}. Then \eqref{equ} becomes
\begin{equation}
  (e^y)_t + (e^y y_{xx})_{xx} = 0, \quad y(\cdot,0)=y_0 \quad\mbox{in }S^1,
    \label{eqy}
\end{equation}
where $y_0=\log u_0$. In order to prove the existence of solutions to this
equation, we semi-discretize \eqref{eqy} in time.
For this, let $T>0$,
and let $0=t_0<t_1<\cdots<t_N=T$ with $t_k=k\tau$ be a partition of $[0,T]$.
Furthermore, let $y_{k-1}\in H^2(S^1)$ with
$\int \exp(y_{k-1})dx=\int u_0dx$ and $\int(\exp(y_{k-1})-y_{k-1})dx
\le\int(u_0-\log u_0)dx$ be given.
Then we solve recursively the elliptic equations
\begin{equation}
  \frac{1}{\tau}(e^{y_k}-e^{y_{k-1}}) + (e^{y_k}(y_k)_{xx})_{xx} = 0
    \quad\mbox{in }S^1.  \label{semi}
\end{equation}

\begin{lemma}\label{ex-disc}
There exists a solution $y_k\in H^2(S^1)$ to \eqref{semi}.
\end{lemma}

\proof Set $z=y_{k-1}$. We consider first for given $\eps>0$ the equation
\begin{equation}
  (e^y y_{xx})_{xx} - \eps y_{xx} + \eps y = \frac{1}{\tau}(e^z-e^y)
    \quad\mbox{in }S^1.  \label{semieps}
\end{equation}
In order to prove the existence of a solution to this approximate problem
we employ the Leray-Schauder theorem. For this, let $w\in H^1(S^1)$
and $\sigma\in[0,1]$ be given, and consider
\begin{equation}
  a(y,\phi) = F(\phi) \quad\mbox{for all }\phi\in H^2(S^1),  \label{semilin}
\end{equation}
where
\begin{eqnarray*}
  a(y,\phi) &=& \int_{S^1}(e^w y_{xx}\phi_{xx}+\eps y_x\phi_x+\eps y\phi)dx,\\
  F(\phi) &=& \frac{\sigma}{\tau}\int_{S^1}(e^z-e^w)\phi dx, \quad
    y,\phi\in H^2(S^1).
\end{eqnarray*}
Clearly, $a(\cdot,\cdot)$ is bilinear, continuous and coercive on $H^2(S^1)$
and $F$ is linear and continuous on $H^2(S^1)$.
(Here we need the additional $\eps$-terms.) Therefore, the Lax-Milgram
lemma provides the existence of a solution $y\in H^2(S^1)$ to \eqref{semilin}.
This defines a fixed-point operator $S:H^1(S^1)\times[0,1]\to H^1(S^1)$,
$(w,\sigma)\mapsto y$. It holds $S(w,0)=0$ for all $w\in H^1(S^1)$. Moreover,
the functional $S$ is continuous and compact
(since the embedding $H^2(S^1)\subset H^1(S^1)$ is compact).
We need to prove a uniform bound for all fixed points of $S(\cdot,\sigma)$.

Let $y$ be a fixed point of $S(\cdot,\sigma)$, i.e.,
$y\in H^2(S^1)$ solves for all $\phi\in H^2(S^1)$
\begin{equation}
  \int_{S^1} (e^y y_{xx}\phi_{xx}+\eps y_x\phi_x+\eps y\phi)dx
    = \frac{\sigma}{\tau}\int_{S^1}(e^z-e^y)\phi dx. \label{semiepsweak}
\end{equation}
Using the test function $\phi=1-e^{-y}$ yields
$$
  \int_{S^1} y_{xx}^2 dx - \int_{S^1} y_{xx}y_x^2dx
    + \eps\int_{S^1}e^{-y}y_x^2dx + \eps\int_{S^1}y(1-e^{-y})dx
    = \frac{\sigma}{\tau}\int_{S^1}(e^z-e^y)(1-e^{-y})dx.
$$
The second term on the left-hand side vanishes since
$y_{xx} y_x^2=(y_x^3)_x/3$.
The third and fourth term on the left-hand side are non-negative.
Furthermore, with the inequality $e^x\ge 1+x$ for all $x\in{\mathbb R}$,
$$
  (e^z-e^y)(1-e^{-y}) \le (e^z-z)-(e^y-y).
$$
We obtain
$$
  \frac{\sigma}{\tau}\int_{S^1}(e^y-y)dx + \int_{S^1}y_{xx}^2 dx
    \le \frac{\sigma}{\tau}\int_{S^1}(e^z-z)dx.
$$
As $z$ is given, this provides a uniform bound for $y_{xx}$ in $L^2(S^1)$.
Moreover, the inequality $e^x-x\ge |x|$ for all $x\in\R$ implies
a (uniform) bound for $y$ in $L^1(S^1)$ and for $\int y dx$.
Now we use the Poincar\'e inequality
$$
  \Big\|u - \int_{S^1}u\frac{dx}L\Big\|_{L^2(S^1)} \le \frac L{2\pi}\|u_x\|_{L^2(S^1)}
    \le \left(\frac L{2\pi}\right)^2\|u_{xx}\|_{L^2(S^1)} \quad\mbox{for all }u\in H^2(S^1).
$$
Recall that $\|u\|_{L^2(S^1)}^2=\int_{S^1}u^2dx/L$.
Then the above estimates provide a (uniform in $\varepsilon$) bound for $y$ and $y_x$
in $L^2(S^1)$ and thus for $y$ in $H^2(S^1)$. This shows that all fixed
points of the operator $S(\cdot,\sigma)$ are uniformly bounded in $H^1(S^1)$.
We notice that we even
obtain a uniform bound for $y$ in $H^2(S^1)$ which is independent of $\eps$.
The Leray-Schauder fixed-point theorem finally ensures the existence of
a fixed point of $S(\cdot,1)$, i.e., of a solution $y\in H^2(S^1)$ to \eqref{semieps}.

It remains to show that the limit $\eps\to 0$ can be performed in
\eqref{semieps} and that the limit function satisfies \eqref{semi}.
Let $y_\eps$ be a solution to \eqref{semieps}. The above estimate shows that
$y_\eps$ is bounded in $H^2(S^1)$ uniformly in $\eps$. Thus there exists
a subsequence (not relabeled) such that, as $\eps\to 0$,
\begin{eqnarray*}
  y_\eps\rightharpoonup y && \mbox{weakly in }H^2(S^1), \\
  y_\eps\to y && \mbox{strongly in }H^1(S^1)\mbox{ and in }L^\infty(S^1).
\end{eqnarray*}
We conclude that $e^{y_\eps}\to e^y$ in $L^2(S^1)$ as $\eps\to 0$.
In particular, $e^{y_\eps}(y_{\eps})_{xx}\rightharpoonup e^y y_{xx}$ weakly
in $L^1(S^1)$.
The limit $\eps\to 0$ in \eqref{semiepsweak} can be performed proving that
$y$ solves \eqref{semi}. Moreover, using the test function $\phi\equiv 1$
in the weak formulation of \eqref{semi} shows that
$\int \exp(y_k)dx=\int\exp(y_{k-1})dx=\int u_0 dx$.
\endproof

For the proof of Theorem \ref{existence} we need further uniform estimates
for the finite sequence $(y^{(N)})$. For this, let $y^{(N)}$ be defined by
$y^{(N)}(x,t)=y_k(x)$ for $x\in S^1$, $t\in(t_{k-1},t_k]$, $1\leq k\leq N$.
Then we have shown in the proof of Lemma \ref{ex-disc}
that there exists a constant $c>0$ depending neither on $\tau$ nor on $N$ such that
\begin{equation}
  \|y^{(N)}\|_{L^2(0,T;H^2(S^1))} + \|y^{(N)}\|_{L^\infty(0,T;L^1(S^1))}
    + \|e^{y^{(N)}}\|_{L^\infty(0,T;L^1(S^1))} \le c. \label{est1}
\end{equation}
To pass to the limit in the approximating equation, we need further compactness estimates on $e^{y^{(N)}}$.
Here we proceed similarly as in \cite{GJT04}.

\begin{lemma}\label{estim2}
The following estimates hold:
\begin{equation}
  \|y^{(N)}\|_{L^{5/2}(0,T;W^{1,\infty}(S^1))}
    + \|e^{y^{(N)}}\|_{L^{5/2}(0,T;W^{1,1}(S^1))} \le c, \label{est2}
\end{equation}
where $c>0$ does not depend on $\tau$ and $N$.
\end{lemma}

\begin{proof} We obtain from the Gagliardo-Nirenberg inequality and
\eqref{est1}:
\begin{eqnarray*}
  \|y^{(N)}\|_{L^{5/2}(0,T;L^\infty(S^1))}
  &\le& c\|y^{(N)}\|_{L^\infty(0,T;L^1(S^1))}^{3/5}
    \|y^{(N)}\|_{L^1(0,T;H^2(S^1))}^{2/5} \le c, \\
  \|y^{(N)}_x\|_{L^{5/2}(0,T;L^\infty(S^1))}
  &\le& c\|y^{(N)}\|_{L^\infty(0,T;L^1(S^1))}^{1/5}
    \|y^{(N)}\|_{L^2(0,T;H^2(S^1))}^{4/5} \le c.
\end{eqnarray*}
This implies the first bound in \eqref{est2}. The second bound follows
from the first one and \eqref{est1}:
\begin{eqnarray*}
  \|e^{y^{(N)}}\|_{L^{5/2}(0,T;W^{1,1}(S^1))}
  &\kern -3pt\le\kern -3pt& c\left(\|e^{y^{(N)}}\|_{L^{5/2}(0,T;L^1(S^1))}
    + \|(e^{y^{(N)}})_x\|_{L^{5/2}(0,T;L^1(S^1))}\right) \\
  &\kern -3pt\le\kern -3pt& c\|e^{y^{(N)}}\|_{L^{5/2}(0,T;L^1(S^1))} + c\|e^{y^{(N)}}\|_{L^\infty(0,T;L^1(S^1))}
    \|y^{(N)}_x\|_{L^{5/2}(0,T;L^\infty(S^1))} \\
  &\kern -3pt\le\kern -3pt& c.
\end{eqnarray*}
The lemma is proved.
\end{proof}

We also need an estimate for the discrete time derivative. We introduce
the shift operator $\sigma_N$ by $(\sigma_N(y^{(N)}))(x,t)=y_{k-1}(x)$
for $x\in S^1$, $t\in(t_{k-1},t_k]$.

\begin{lemma}\label{estim3}
The following estimate holds:
\begin{equation}
  \|e^{y^{(N)}}-e^{\sigma_N(y^{(N)})}\|_{L^1(0,T;H^{-2}(0,1))}
    \le c\tau,  \label{est3}
\end{equation}
where $c>0$ does not depend on $\tau$ and $N$.
\end{lemma}

\begin{proof}
{}From \eqref{semi} and H\"older's inequality we obtain
\begin{eqnarray*}
  &\frac{1}{\tau}\|e^{y^{(N)}}-e^{\sigma_N(y^{(N)})}\|_{L^1(0,T;H^{-2}(S^1))}
    \le \|e^{y^{(N)}}y^{(N)}_{xx}\|_{L^1(0,T;L^2(S^1))}\\
  &\le \|e^{y^{(N)}}\|_{L^2(0,T;L^\infty(S^1))}
    \|y^{(N)}_{xx}\|_{L^2(0,T;L^2(S^1))},
\end{eqnarray*}
and the right-hand side is uniformly bounded by \eqref{est1} and
\eqref{est2} since $W^{1,1}(0,1)\hookrightarrow L^\infty(0,1)$.
\end{proof}

Now we are able to prove Theorem \ref{existence}, i.e.\ to perform the
limit $\tau\to 0$ in \eqref{semi}.
{}From estimate \eqref{est1} the existence of a subsequence
of $y^{(N)}$ (not relabeled) follows
such that, as $N\to\infty$ or, equivalently, $\tau\to 0$,
\begin{equation}
  y^{(N)}\rightharpoonup y \quad\mbox{weakly in }L^2(0,T;H^2(S^1)).
    \label{conv1}
\end{equation}
Since the embedding $W^{1,1}(S^1)\subset L^1(S^1)$ is compact
it follows from the second bound in \eqref{est2} and from \eqref{est3}
by an application of Aubin's lemma \cite[Thm.~5]{Sim87} that, up to the extraction of
a subsequence, $e^{y^{(N)}}\to g$ strongly in $L^1(0,T;L^1(S^1))$ and hence
also in $L^1(0,T;H^{-2}(S^1))$.

We claim that $g=e^y$. For this, we observe that, by \eqref{est1},
\begin{eqnarray}
  \|e^{y^{(N)}}-g\|_{L^2(0,T;H^{-2}(S^1))}^2
  &\le& \|e^{y^{(N)}}-g\|_{L^\infty(0,T;H^{-2}(S^1))}
    \|e^{y^{(N)}}-g\|_{L^1(0,T;H^{-2}(S^1))} \nonumber \\
  &\le& c\left(\|e^{y^{(N)}}\|_{L^\infty(0,T;L^1(S^1))}
    + \|g\|_{L^\infty(0,T;L^1(S^1))}\right) \nonumber \\
  && {}\qquad\times \|e^{y^{(N)}}-g\|_{L^1(0,T;H^{-2}(S^1))} \nonumber \\
  &\le& c\|e^{y^{(N)}}-g\|_{L^1(0,T;H^{-2}(S^1))} \to 0
    \quad\mbox{as }N\to\infty.  \nonumber
\end{eqnarray}
Now let $z$ be a smooth function. Since $e^{y^{(N)}}\to g$ strongly in
$L^2(0,T;H^{-2}(S^1))$ and $y^{(N)}\rightharpoonup y$ weakly
in $L^2(0,T;H^2(S^1))$, we can pass to the limit $N\to\infty$ in
$$
  0\le \int_0^T\langle e^{y^{(N)}}-e^z,y^{(N)}-z\rangle_{H^{-2},H^2} dt
$$
to obtain the inequality
$$
  0\le\int_0^T\int_{S^1}(g-e^z)(y-z)dxdt.
$$
The monotonicity of $x\mapsto e^x$ finally yields $g=e^y$.

In particular, $e^{y^{(N)}}\to e^y$ strongly in $L^1(0,T;L^1(S^1))$.
The second uniform bound in~\eqref{est2} implies that, up to the possible extraction of a subsequence again, $e^{y^{(N)}}\rightharpoonup e^y$ weakly* in
$L^{5/2}(0,T;$ $L^\infty(S^1))$. Thus, Lebesgue's convergence theorem
gives
\begin{equation}
  e^{y^{(N)}} \to e^y \quad\mbox{strongly in }L^2(0,T;L^2(S^1)).
    \label{conv2}
\end{equation}
Furthermore, the uniform estimate \eqref{est3} implies, for a
subsequence,
\begin{equation}
  \frac{1}{\tau}\left(e^{y^{(N)}}-e^{\sigma_N(y^{(N)})}\right)
    \rightharpoonup (e^y)_t \quad\mbox{weakly in }L^1(0,T;H^{-2}(S^1)).
    \label{conv3}
\end{equation}
We can pass to the limit $\tau\to 0$ in \eqref{semi}, using the
convergence results \eqref{conv1}-\eqref{conv3}, which concludes the proof of Theorem~\ref{existence}.\endproof


\section{Uniqueness of solutions}\label{sec-un}

To get a uniqueness result, we need an additional
regularity assumption.

\begin{theorem}\label{unique}
Let $u_1$, $u_2$ be two weak solutions to \eqref{equ} in the sense of
Theorem \ref{existence} with the same initial data
such that $u_1,u_2\in C^0([0,T];L^1(S^1))$ and
$\sqrt{u_1/u_2}$, $\sqrt{u_2/u_1}\in
L^2(0,T;$ $H^2(S^2))$ for some $T>0$. Then $u_1=u_2$ in $S^1\times(0,T)$.
\end{theorem}

Bleher et al.\ have showed the uniqueness of solutions to \eqref{equ}
in the class of mild solutions, i.e.\ $C^0([0,T];H^1(S^1))$, which are
{\em positive}. We allow for the more general class of {\em non-negative} solutions satisfying the above regularity assumptions.

\proof We use a similar idea as in \cite{JuPi00}. Employing
the test function $1-\sqrt{u_2/u_1}$ in equation~\eqref{equ} for $u_1$ and the test function $\sqrt{u_1/u_2}-1$ in equation~\eqref{equ} for $u_2$ and taking the difference of both equations yields
\begin{eqnarray*}
  \lefteqn{\int_0^t\left\langle (u_1)_t,1-\sqrt{\frac{u_2}{u_1}}
    \right\rangle_{H^{-2},H^2} dt
    - \int_0^t\left\langle (u_2)_t,\sqrt{\frac{u_1}{u_2}}-1
    \right\rangle_{H^{-2},H^2} dt} \\
  &=& \int_0^t\left\langle(u_1(\log u_1)_{xx})_{xx},\sqrt{\frac{u_2}{u_1}}
    \right\rangle_{H^{-2},H^2} dt + \int_0^t\left\langle(u_2(\log u_2)_{xx})_{xx},\sqrt{\frac{u_1}{u_2}}
    \right\rangle_{H^{-2},H^2} dt \\
  &=& I_1 + I_2.
\end{eqnarray*}
The left-hand side can be {\em formally} written as
\begin{eqnarray*}
  \lefteqn{\int_0^t\left\langle (u_1)_t,1-\sqrt{\frac{u_2}{u_1}}
    \right\rangle_{H^{-2},H^2} dt
    - \int_0^t\left\langle (u_2)_t,\sqrt{\frac{u_1}{u_2}}-1
    \right\rangle_{H^{-2},H^2} dt} \\
  &=& 2\int_0^t\int_{S^1}\left[(\sqrt{u_1})_t(\sqrt{u_1}-\sqrt{u_2})
    - (\sqrt{u_2})_t(\sqrt{u_1}-\sqrt{u_2})\right]dxdt \\
  &=& \int_{S^1}\left(\sqrt{u_1(t)}-\sqrt{u_2(t)}\right)^2 dx.
\end{eqnarray*}
As the first and the last equation hold rigorously, it is possible
to make the computation rigorous by approximating $u_1$ and $u_2$
by suitable smooth functions and then passing to the limit in the
first and the last equation by a standard procedure.

We claim now that $I_1+I_2$ is non-positive. For this we compute
{\em formally} as follows.
\begin{eqnarray*}
  I_1 &=& 2\int_0^t\left\langle (\sqrt{u_1})_{xxxx}-\frac{1}{\sqrt{u_1}}
    \left|(\sqrt{u_1})_{xx}\right|^2, \sqrt{u_2}\right\rangle_{H^{-2},H^2}dt\\
  &=& -2\int_0^t\int_{S^1}\left[-(\sqrt{u_1})_{xx}(\sqrt{u_2})_{xx}
    + \left|(\sqrt{u_1})_{xx}\right|^2\sqrt{\frac{u_2}{u_1}}\,\right]dxdt.
\end{eqnarray*}
A similar result can be obtained for $I_2$. Thus
$$
  I_1+I_2 = -2\int_0^t\int_{S^1}\left|\sqrt[4]{\frac{u_2}{u_1}}
    (\sqrt{u_1})_{xx} - \sqrt[4]{\frac{u_1}{u_2}}
    (\sqrt{u_2})_{xx}\right|^2 \le 0.
$$
This calculation can be made rigorous again by an approximation argument.
We conclude that
$$
  \int_{S^1}\left|\sqrt{u_1(t)}-\sqrt{u_2(t)}\right|^2 dx \le 0,
$$
which gives $u_1(t)=u_2(t)$ in $S^1$ for all $t\le T$.
\endproof


\section{Optimal logarithmic Sobolev inequality on $S^1$}\label{sec-so}

The main goal of this section is the proof of a logarithmic
Sobolev inequality for periodic functions.
The following theorem is due to Weissler and Rothaus
(see \cite{EY87,Ro80,Wei80}). We give a simple proof using the
entropy--entropy production method.
Recall that $S^1$ is parametrized by $0\le x\le L$.

\begin{theorem}\label{logsob}
Let $\Hh_1=\{u\in H^1(S^1):u_x\not\equiv 0\mbox{\rm\ a.e.}\}$ and $\|u\|_{L^2(S^1)}^2=
\int_{S^1}u^2 dx/L$. Then
\begin{equation}
  \inf_{u\in \Hh_1}\frac{\int_{S^1}u_x^2 dx}{\int_{S^1}u^2
    \log(u^2/\|u\|_{L^2(S^1)}^2)dx} = \frac{2\pi^2}{L^2}.  \label{optimal}
\end{equation}
\end{theorem}

\noindent We recall that the optimal constant in the usual Poincar\'e inequality
is $L/2\pi$, i.e.
\begin{equation}
   \inf_{v\in \Hh_1}\frac{\int_{S^1}v_x^2 dx}{\int_{S^1}(v-\bar v)^2 dx}
    = \frac{4\pi^2}{L^2},  \label{poincare}
\end{equation}
where $\bar v=\int_{S^1}v dx/L$.

\proof Let $I$ denote the value of the infimum in \eqref{optimal}.
Let $u\in \Hh_1$ and define $v$ by setting $u=1+\eps(v-\bar v)$. Then, if
we can prove that
\begin{equation}
  I \le \frac12\inf_{v\in \Hh_1}\frac{\int_{S^1}v_x^2 dx}{\int_{S^1}
    (v-\bar v)^2 dx}, \label{aux1}
\end{equation}
we obtain the upper bound $I\le 2\pi^2/L^2$ from \eqref{poincare}.
Without loss of generality, we may replace $v-\bar v$ by $v$ such that
$\int_{S^1}vdx = 0$. Then $u^2 = 1+2\eps v+\eps^2 v^2$ and the
expansion $\log(1+x)=x+x^2/2+O(x^3)$ for $x\to 0$ yield for $\eps\to 0$
\begin{eqnarray*}
  \int_{S^1}u^2\log(u^2) dx &=& \int_{S^1}(1+2\eps v+\eps^2 v^2)
    \log(1+2\eps v+\eps^2 v^2)dx \\
  &=& 3\eps^2\int_{S^1}v^2 dx + O(\eps^3), \\
  \int_{S^1}u^2dx\log\left(\frac{1}{L}\int_{S^1}u^2 dx\right)
  &=& \int_{S^1}(1+\eps^2v^2)dx\log\left(\frac{1}{L}
    \int_{S^1}(1+\eps^2v^2)dx\right) \\
  &=& \eps^2\int_{S^1}v^2 dx + O(\eps^4).
\end{eqnarray*}
Taking the difference of the two expansions gives
$$
  \int_{S^1}u^2\log\left(\frac{u^2}{\int_{S^1}u^2 dx/L}\right)dx
   = 2\eps^2\int_{S^1}v^2 dx + O(\eps^3).
$$
Therefore, using $\int_{S^1}u_x^2 dx = \eps^2\int_{S^1}v_x^2 dx$,
$$
  \frac{\int_{S^1}u_x^2 dx}{\int_{S^1}u^2
    \log(u^2/\|u\|_{L^2(S^1)}^2)dx} = \frac12\,\frac{\int_{S^1}v_x^2dx}{\int_{S^1}v^2dx}
    + O(\eps).
$$
In the limit $\eps\to 0$ we obtain \eqref{aux1}.

In order to prove the lower bound for the infimum we use the
entropy--entropy production method. For this we consider the heat equation
$$
  v_t = v_{xx}\quad\mbox{in }S^1\times(0,\infty), \quad
    v(\cdot,0)=u^2\quad\mbox{in }S^1
$$
for some function $u\in H^1(S^1)$. We assume for simplicity that
$\|u\|_{L^2(S^1)}^2=\int_{S^1}u^2dx/L=1$. Then
$$
  \frac{d}{dt}\int_{S^1}v\log v dx = -4\int_{S^1}w_x^2 dx,
$$
where the function $w:=\sqrt{v}$ solves the equation $w_t=w_{xx}+w_x^2/w$.
Now, the time derivative of
$$
  f(t)=\int_{S^1}w_x^2 dx - \frac{2\pi^2}{L^2}\int_{S^1}w^2\log(w^2)dx
$$
equals
$$
  f'(t) = -2\int_{S^1}\left(w_{xx}^2+\frac{w_x^4}{3w^2}
    -\frac{4\pi^2}{L^2} w_x^2\right)dx
  \le -\frac23\int_{S^1}\frac{w_x^4}{w^2}dx \le 0,
$$
where we have used the Poincar\'e inequality
\begin{equation}
  \int_{S^1}w_x^2 dx \le \frac{L^2}{4\pi^2}\int_{S^1}w_{xx}^2 dx.
    \label{star}
\end{equation}
This shows that $f(t)$ is non-increasing and moreover,
for any $u\in H^1(S^1)$,
$$
  \int_{S^1}u_x^2dx - \frac{2\pi^2}{L^2}\int_{S^1}u^2\log(u^2/\|u\|_{L^2(S^1)}^2)dx
    = f(0) \ge f(t).
$$
As the solution $v(\cdot,t)$ of the above heat equation and hence
$w(\cdot,t)$ converges to zero in appropriate Sobolev norms as $t\to+\infty$,
we conclude that $f(t)\to 0$ as $t\to+\infty$. This implies $I\ge2\pi^2/L^2$.
\endproof

\begin{remark}\rm
Similar results as in Theorem \ref{logsob} can be obtained for the so-called
convex Sobolev inequalities. Let $\sigma(v)=(v^p-\bar v^p)/(p-1)$,
where $\bar v=\int_{S^1}vdx/L$ for $1<p\le 2$. We claim that
$$
  \inf_{v\in \Hh_1}\frac{\int_{S^1}\sigma''(v)v_x^2 dx}{\int_{S^1}\sigma(v)dx}
    = \frac{8\pi^2}{L^2}.
$$
As in the logarithmic case, the lower bound is achieved by an expansion
around~$1$ and the usual Poincar\'e inequality. On the other hand,
let $v$ be a solution of the heat equation. Then
$$
\frac d{dt}\int_{S^1} \sigma(v)dx=-\frac 4p\int_{S^1} w_x^2dx
$$
where $w=v^{p/2}$ solves
\begin{equation}\label{Eqn:Sqrt}
w_t=w_{xx}+\left(\frac2p-1\right)\frac{w_x^2}w,
\end{equation}
and, using \eqref{star},
\begin{eqnarray*}
  \frac d{dt}\int_{S^1} \left(w_x^2-\frac{2\pi^2p}{L^2}\sigma(v)\right)dx
  &=& -2\int_{S^1} \left(w_{xx}^2-\frac{4\pi^2}{L^2}w_x^2
    +\left(\frac2p-1\right)\frac{w_x^4}{3w^2}\right)dx \\
  &\leq& -\frac 23\left(\frac2p-1\right)\int_{S^1}\frac{w_x^4}{w^2}dx \le 0.
\end{eqnarray*}
This proves the upper bound
$$
  \frac p4\int_{S^1}\sigma''(v)v_x^2 dx = \int_{S^1} w_x^2dx
    \geq \frac{2\pi^2p}{L^2}\int_{S^1}\sigma(v) dx.
$$

With the notation $v=u^{2/p}$ this result takes the more familiar form
\begin{equation}\label{convex}
  \frac{1}{p-1}\left[\int_{S^1}u^2dx - L\left(\frac 1L\int_{S^1}u^{2/p}dx
    \right)^p\right] \le \frac{L^2}{2\pi^2p}\int_{S^1}u_x^2dx
    \quad\mbox{for all } u\in H^1(S^1).
\end{equation}
The logarithmic case corresponds to the limit $p\to 1$ whereas the
case $p=2$ gives the usual Poincar\'e inequality.

We may notice that the method gives more than what is stated in Theorem \ref{logsob} since there is an integral remainder term. Namely, for any $p\in[1,2]$, for any $v\in H^1(S^1)$, we have
$$
  \frac p4\int_{S^1}\sigma''(v)v_x^2 dx +{\cal R}[v]
    \geq \frac{2\pi^2p}{L^2}\int_{S^1}\sigma(v) dx
$$
with
$$
{\cal R}[v]=2\int_0^\infty\int_{S^1} \left(w_{xx}^2-\frac{4\pi^2}{L^2}w_x^2
    +\left(\frac2p-1\right)\frac{w_x^4}{3w^2}\right)dx\,dt,
$$
where $w=w(x,t)$ is the solution to \eqn{Eqn:Sqrt} with initial datum $u_0^{p/2}$. Inequality \eqn{convex} can also be improved with an integral remainder term for any $p\in[1,2]$, where in the limit case $p=1$, one has to take $\sigma(v)=v\log(v/\bar v)$. As a consequence, the only optimal functions in \eqn{optimal} or in \eqn{convex} are the constants.
\end{remark}

\begin{corollary}\label{logsob2}
Let $n\in\N$, $n>0$ and let $\Hh_n=\{u\in H^n(S^1):u_x\not\equiv 0\mbox{\rm\
a.e.}\}$. Then
\begin{equation}
  \inf_{u\in \Hh_n}\frac{\int_{S^1}\left|u^{(n)}\right|^2 dx}{\int_{S^1}u^2
    \log(u^2/\|u\|_{L^2(S^1)}^2)dx} = \frac{1}{2}\left(\frac{2\pi}{L}\right)^{2n}.
    \label{optimal2}
\end{equation}
\end{corollary}

\proof We obtain a lower bound by applying successively Theorem \ref{logsob}
and the Poincar\'e inequality:
\[
\int_{S^1}u^2\log\left(\frac{u^2}{\|u\|_{L^2(S^1)}^2}\right)dx
\leq \frac{L^2}{2\pi^2}\int_{S^1}u_x^2dx\leq
2\left(\frac{L}{2\pi}\right)^{2n}\int_{S^1}\left|u^{(n)}\right|^2 dx
\]

The upper bound is achieved as in the proof of Theorem \ref{logsob}
by expanding the quotient for $u= 1+\eps v$ with $\int_{S^1}v dx=0$
in powers of $\eps$,
\[
\frac{\int_{S^1}\left|u^{(n)}\right|^2 dx}{\int_{S^1}u^2 \log(u^2/\|u\|_{L^2(S^1)}^2)dx}
=\frac 12\,\frac{\int_{S^1}\left|v^{(n)}\right|^2 dx}{\int_{S^1}v^2dx} +O(\eps),
\]
and using the Poincar\'e inequality
\[
\inf_{u\in \Hh_n}\frac{\int_{S^1}\left|v^{(n)}\right|^2 dx}{\int_{S^1}\left|v-\bar v\right|^2dx}
=\left(\frac{2\pi}{L}\right)^{2n}.
\]
The best constant $\omega=\left({2\pi}/{L}\right)^{2n}$ in such an
inequality is easily recovered by looking
for the smallest positive value of $\omega$ for which there exists a  nontrivial periodic solution of $(-1)^nv^{(2n)}+\omega v=0$.
\endproof


\section{Exponential time decay of the solutions}\label{sec-lo}

We show the exponential time decay of the solutions of \eqref{equ}.
Our main result is contained in the following theorem.

\begin{theorem}\label{longtime}
Assume that $u_0$ is a nonnegative measurable function such that $\int_{S^1}
(u_0-\log u_0) dx$ and $\int_{S^1}u_0\log u_0 dx$ are finite.
Let $u$ be the weak solution of \eqref{equ} constructed in Theorem
\ref{existence} and set $\bar u=\int_{S^1}u_0(x)
dx/L$. Then
$$
  \int_{S^1}u(\cdot,t)\log\left(\frac{u(\cdot,t)}{\bar u}\right)dx \le
    e^{-Mt}\int_{S^1}u_0\log\left(\frac{u_0}{\bar u}\right)dx,
$$
where
$$
  M = \frac{32\pi^4}{L^4}.
$$
\end{theorem}

\proof Since we do not have enough regularity of the solutions to
\eqref{equ} we need to regularize the equation first. For this we consider the semi-discrete
problem
\begin{equation}\label{semi2}
  \frac{1}{\tau}(u_k-u_{k-1}) + (u_k(\log u_k)_{xx})_{xx} = 0\quad\mbox{in }
    S^1
\end{equation}
as in the proof of Theorem \ref{existence}. The solution $u_k\in H^2(S^1)$
of this problem for given $u_{k-1}$ is strictly positive and we can use
$\log u_k$ as a test function in the weak formulation of \eqref{semi2}.
In order to simplify the presentation we set $u:=u_k$ and $z:=u_{k-1}$.
Then we obtain as in~\cite{JuTo03}
\begin{equation}\label{semi2a}
  \frac{1}{\tau}\int_{S^1}(u\log u - z\log z)dx + \int_{S^1}
    u\left|(\log u)_{xx}\right|^2 dx \le 0.
\end{equation}
{}From integration by parts it follows
$$
  \int_{S^1}\frac{u_x^2 u_{xx}}{u^2}dx = \frac23\int_{S^1}
    \frac{u_x^4}{u^3}dx.
$$
This identity gives
\begin{eqnarray*}
  \int_{S^1} u\left|(\log u)_{xx}\right|^2 dx
  &=&\! \int_{S^1}\!\left(\frac{u_{xx}^2}{u}+\frac{u_x^4}{u^3}
    -2\frac{u_{xx}u_x^2}{u^2}\right)dx
    = \int_{S^1}\left(\frac{u_{xx}^2}{u}-\frac13\,\frac{u_x^4}{u^3}\right)dx \\
  &=&\! \int_{S^1}\!\left(\frac{u_{xx}^2}{u}+\frac13\,\frac{u_x^4}{u^3}
    -\frac{u_{xx}u_x^2}{u^2}\right)dx
    = 4\int_{S^1}|(\sqrt{u})_{xx}|^2 dx + \frac{1}{12}
    \int_{S^1}\frac{u_x^4}{u^3}dx.
\end{eqnarray*}
Thus, \eqref{semi2a} becomes
\begin{equation}\label{semi3}
  \frac{1}{\tau}\int_{S^1}\Big(u\log\left(\frac{u}{\bar u}\right) - z\log\left(\frac{z}{\bar u}\right)
    \Big)dx + 4\int_{S^1}\left|(\sqrt{u})_{xx}\right|^2dx \le 0.
\end{equation}
Now we use Corollary~\ref{logsob2} with $n=2$:
$$
  \int_{S^1}u\log\left(\frac{u}{\bar u}\right)dx \le \frac{L^4}{8\pi^4}
   \int_{S^1}\left|(\sqrt{u})_{xx}\right|^2dx.
$$
{}From this inequality and \eqref{semi3} we conclude
$$
  \frac{1}{\tau}\int_{S^1}\Big(u\log\left(\frac{u}{\bar u}\right) - z\log\left(\frac{z}{\bar u}\right)
    \Big)dx + \frac{32\pi^4}{L^4}\int_{S^1}u\log\left(\frac{u}{\bar u}\right)dx
    \le 0.
$$
This is a difference inequality for the sequence
$$
  E_k := \int_{S^1}u_k\log\left(\frac{u_k}{\bar u}\right)dx,
$$
yielding
$$
  (1+\tau M)E_k\le E_{k-1}\quad\mbox{or}\quad E_k\le E_0(1+\tau M)^{-k},
$$
where $M$ is as in the statement of the theorem. For $t\in((k-1)\tau,k\tau]$
we obtain further
$$
  E_k \le E_0(1+\tau M)^{-t/\tau}.
$$
Now the proof as exactly as in \cite{JuTo03}. Indeed, the functions $u_k(x)$
converge a.e.\ to $u(x,t)$ and $(1+\tau M)^{-t/\tau}\to e^{-Mt}$ as
$\tau\to 0$. This implies the assertion.
\endproof

\begin{remark}\rm
The decay rate $M$ is not optimal since in the estimate \eqref{semi3} we have neglected
the term
$\frac{1}{12}\int(u_x^4/u^3)dx$ .
\end{remark}


\end{document}